\documentclass[12pt]{article}
\usepackage{amsmath,amsfonts}

\newtheorem{Theorem}{Theorem}[section]

\newtheorem{Proposition}[Theorem]{Proposition}
\newtheorem{Lemma}[Theorem]{Lemma}
\newtheorem{Corollary}[Theorem]{Corollary}
\newtheorem{Remark}[Theorem]{Remark}

\newtheorem{Hypothesis}[Theorem]{Hypothesis}
\makeatletter
\newif\ifmsbmloaded@

\@addtoreset{equation}{section}

\makeatother

\def\hph#1{\hphantom{#1}}

\def\R{\mathbb R}
\def\N{\mathbb N}

\def\E{\mathbb E}
\def\P{\mathbb P}
\def\ds{\displaystyle}

\begin{document}

\begin{center} {\Large \bf Coupling for some partial   differential equations driven by white noise}
\vspace{5mm}

 Giuseppe Da Prato, Scuola Normale Superiore di Pisa,

 Piazza dei Cavalieri 7, 56126, Pisa, Italy\vspace{2mm}

 Arnaud Debussche, Ecole Normale Sup\'erieure de Cachan,
  antenne de Bretagne,

 Campus de Ker Lann, 35170 Bruz, France.\vspace{2mm}

 Luciano Tubaro, Department of Mathematics, University of Trento,

 Via Sommarive  14, 38050 Povo, Italy.

\end{center}\vspace{5mm}

\begin{abstract} We prove, using coupling arguments,   exponential convergence to equilibrium for reaction--diffusion and
Burgers  equations driven by space-time white
noise. We use a coupling by reflection.

\end{abstract}



\section{Introduction} 

We are concerned with a  stochastic differential  equation in a separable Hilbert space $H$, 
with inner product 
$(\cdot,\cdot  )$ and norm $|\cdot|$,
\begin{equation}
\label{e1.1}
dX=(AX+b(X))dt+dW(t),\quad X(0)=x\in H,
\end{equation}
where $A\colon D(A)\subset H\to H$ is linear,   $b\colon D(b)\subset H\to H$ is nonlinear
and $W$ is a cylindrical Wiener process defined in some probability space $(\Omega, \mathcal F,\P)$ in $H$.
Concerning $A$ we shall assume that
\begin{Hypothesis}
\label{h1.1}
\begin{enumerate}
\item[]
\item[(i)] $A\colon D(A)\subset H\to H$ is the infinitesimal
generator of a strongly continuous semigroup $e^{tA}$.

\item[(ii)]  For any $t> 0$  the linear operator $Q_{t},$ defined as  
\begin{equation}
\label{e1.2}
Q_{t}x=\int_{0}^{t}e^{sA}e^{sA^{*}}xds,\;\;x\in H,\;t\ge 0,
\end{equation}
is of trace class.
\end{enumerate}
\end{Hypothesis} 

We shall consider situations where \eqref{e1.1} has a unique mild solution 
 $X(\cdot,x)$, that is a mean square adapted process
such that
\begin{equation}
\label{e1.3}
X(t,x)=e^{tA}x+\int_0^te^{(t-s)A}b(X(s,x))ds+z(t),\quad\P\mbox{\rm--a.s.},
\end{equation} 
where $z(t)$ is the stochastic convolution
\begin{equation}
\label{e1.4}
z(t)=\int_0^te^{(t-s)A}dW(s).
\end{equation} 
It is well known that, thanks to Hypothesis \ref{h1.1}, for each $t> 0$, $z(t)$ is a 
Gaussian random variable with mean $0$ and
covariance operator
$Q_t.$

We will also assume that the solution has continuous trajectories. More precisely, we assume
\begin{equation}
\label{e1.4bis}
X(\cdot,x)\in L^2(\Omega;C([0,T];H)),\mbox{ for any } x\in H.
\end{equation}

In this paper we want to study the exponential convergence to equilibrium of the transition semigroup
$$
P_t\varphi(x)=\E[\varphi(X(t,x))],\quad x\in H,\;t\ge 0,
$$
where $  \varphi\colon H\to \R$. We wish to use coupling arguments.  It is well known that exponential convergence to equilibrium
implies the uniqueness of invariant measure.

It seems that the first paper using a coupling method to prove uniqueness of the invariant 
measure and mixing
property for a stochastic partial differential equation is  \cite{mueller}. There, an equation with globally Lipschitz coefficients is considered, 
some of them are also assumed to be monotone.

Coupling argument have also been used recently  to prove ergodicity and exponential convergence to 
equilibrium in the context of the Navier-Stokes equation driven by very degenerate noises
(see \cite{KS1}, \cite{ks}, \cite{mattingly1}, \cite{mattingly2}). The method has also 
been studied in   \cite{Hairer1} for reaction--diffusion equations and  in  \cite{odasso}
for Ginzburg--Landau equations.

Our interest here is different since we are interested in space-time white noises as in 
 \cite{mueller} but without the strong restrictions on the coefficients. Ergodicity
is well known in the cases considered here. It can be proved 
 by the Doob theorem (see \cite{DPZ2}). Indeed, 
since the noise is non degenerate, it is not diffiucult to prove that the transition semigroup
is strong Feller and irreducible. However, this argument does not imply exponential 
convergence to equilibrium and we think that it is important to study this question.

In this paper we shall follow the construction of couplings  introduced in \cite{LR} 
(see also \cite{CL}) to treat both reaction diffusion
and Burgers equations driven by white noise and  obtain exponential 
convergence to equilibrium. 

Note that  exponential convergence to equilibrium  for reaction--diffusion equations  is well 
known. Anyway, we have chosen to treat this example because we think that the method
presented here provides a very simple proof. Moreover, we recover the spectral gap property 
obtained in \cite{DPDG} by a totally different
- and simpler - argument.

In the case of  the  Burgers equation driven by space-time white noise, it seems that 
our result is new.  

The coupling method based on Girsanov transform introduced  in 
\cite{KS1}, \cite{ks}, \cite{mattingly1}, \cite{mattingly2} can easily be used 
if the noise is nuclear. It is also possible  
it could   be extended to our case. However, the extension is not straightforward and 
the method used here seems to be simpler. 
Moreover, it is not clear that, in the case of the reaction-diffusion
equation, it is possible to prove the spectral gap property with this method.

Next section is devoted to describing the  construction of the coupling  used here, we follow \cite{CL}.
Note that the coupling is constructed as the solution of a stochastic differential 
equation with  discontinuous coefficients. In \cite{CL}, the existence of the coupling is 
straightforward. It is easy to see that the corresponding martingale problem has a solution.
This argument is difficult in infinite dimension and we have preferred to prove directly
the existence of a strong solution.
Section 3  is devoted to application to reaction--diffusion equations and 
section 4 to the Burgers equation driven by white noise.

We finally remark that our method extends to other equations such as reaction-diffusion equations
or the stochastic Navier-Stokes equation in space dimension two with non degenerate noise. We 
have chosen
to restrict our presentation to these two examples for clarity of the presentation.

\section{Construction of the coupling}

We shall consider the following  system of stochastic differential equations:
\begin{equation}
\label{e1.5}
\left\{\begin{array}{lll}
dX_1=(AX_1+b(X_1))dt+\frac{1}{\sqrt{2}}\;dW_1+\frac{1}{\sqrt{2}}\;\left(1-2\frac{(X_1-X_2)\otimes(X_1-X_2)}{|X_1-X_2|^2}\right)dW_2\\
\\
dX_2=(AX_2+b(X_2))dt+\frac{1}{\sqrt{2}}\;\left(1-2\frac{(X_1-X_2)\otimes(X_1-X_2)}{|X_1-X_2|^2}\right)dW_1+\frac{1}{\sqrt{2}}\;dW_2\\
\\
X_1(0)=x_1,\quad X_2(0)=x_2,
\end{array}\right. 
\end{equation}
where $W_1,W_2$ are independent cylindrical Wiener processes.
This corresponds to a coupling with reflection, see \cite{CL}. Equation 
\eqref{e1.5} is associated to  the Kolmogorov operator
   in  $H\times H$ defined by
$$
\begin{array}{lll}
\mathcal K \;\Phi(x_1,x_2)\\
\\
=\ds{\frac12\;\mbox{\rm Tr}\;\left[ \left(\begin{array}{cc}
1&{\ds 1-2\frac{(x_1-x_2)\otimes(x_1-x_2)}{|x_1-x_2|^2}}\\
\ds{1-2\frac{(x_1-x_2)\otimes(x_1-x_2)}{|x_1-x_2|^2}}&1
\end{array}   \right)
 D^2 \Phi(x_1,x_2)\right]}\\
\\
+\left(\left(\begin{array}{c}
Ax_1+b(x_1)\\
Ax_2+b(x_2)
\end{array}   \right),D \Phi(x_1,x_2)  \right),
\end{array}
$$
or, equivalently,
\begin{equation}
\label{e1.6}
\begin{array}{lll}
\mathcal K \;\Phi
&=&\ds{\frac12\;\mbox{\rm Tr}\;\left[\Phi_{x_1x_1}+\left(2-4\frac{(x_1-x_2)\otimes(x_1-x_2)}{|x_1-x_2|^2}\right)
\Phi_{x_1x_2}+\Phi_{x_2x_2}
\right]}\\
\\
&&+( Ax_1+b(x_1),\Phi_{x_1}   )+( Ax_2+b(x_2),\Phi_{x_2}   ).
\end{array}
\end{equation}
The following formula will be useful in the sequel.
\begin{Lemma}
\label{l1.1}
Let $f\colon\R\mapsto\R$ be a ${\mathcal C}^2$ function and let $\Phi$ 
be defined by $\Phi(x_1,x_2)=f(|x_1-x_2|),\; x_1,x_2\in H,\;f\in C^2(\R)$, then 

\begin{equation}
\label{e1.13}
\begin{array}{lll}
\mathcal K \;\Phi(x_1,x_2)
&=&\ds{2\;f''(|x_1-x_2|)+ \frac{f'(|x_1-x_2|)}{|x_1-x_2|}}\\
\\
&&\times ( A(x_1-x_2)+ b(x_1)-b(x_2),x_1-x_2 ).
\end{array}
\end{equation}

\end{Lemma}
{\bf Proof:}
We have
$$
\Phi_{x_1}=-\Phi_{x_2},\quad \Phi_{x_1x_1}=\Phi_{x_2x_2}=-\Phi_{x_1x_2},
$$
$$
\Phi_{x_1}(x_1,x_2)=f'(|x-y|)\;\frac{x_1-x_2}{|x_1-x_2|}, 
$$
and
$$
\begin{array}{lll}
\Phi_{x_1x_1}(x_1,x_2)&=&\ds{f'(|x_1-x_2|)\;\frac{ |x_1-x_2|^2-(x_1-x_2)\otimes(x_1-x_2)}{|x_1-x_2|^3}}\\
\\
&&\ds{+f''(|x_1-x_2|)\; \frac{(x_1-x_2)\otimes(x_1-x_2)}{|x_1-x_2|^2},}
\end{array}
$$
The result follows.
\hfill$\Box$\bigskip

We will use functions $f$ 
such that, for a suitable positive constant $\kappa$, we have  
\begin{equation}
\label{e1.12}
\mathcal K\Phi(x_1,x_2)\le -\kappa,\quad \mbox{\rm for all}\;x_1,x_2\in H.
\end{equation}
Thanks to Lemma \ref{l1.1}, we have to solve the following  basic inequality,
 in the unknown $f$ (notice that $f$ has to be nonnegative), 
\begin{equation}
\label{e1.14}
 \ds{2\;f''(|x_1-x_2|)}
\ds{+ \frac{f'(|x_1-x_2|)}{|x_1-x_2|}\; ( A(x_1-x_2)+ b(x_1)-b(x_2),x_1-x_2   )}
\le -\kappa.
\end{equation}

We now study problem \eqref{e1.5}.
In our applications, it will  be easy to verify that, 
for any $\varepsilon >0$, it has a unique mild solution 
$X(t,x_1,x_2)=(X_1(t,x_1,x_2),X_2(t,x_1,x_2))$ on the random interval
$[0,\tau_{x_1,x_2}^\varepsilon]$ where 
$\tau_{x_1,x_2}^\varepsilon=\inf\{t\in [0,T]\,|\, |X_1(t)-X_2(t)|\le \varepsilon\}$.

Clearly $\tau_{x_1,x_2}^\varepsilon$ is increasing as $\varepsilon\to 0$ so that we can 
define $\tau_{x_1,x_2}=\lim_{\varepsilon\to 0}\tau_{x_1,x_2}^\varepsilon$ and 
get a unique mild solution 
$X(t,x_1,x_2)=(X_1(t,x_1,x_2),X_2(t,x_1,x_2))$ on the random interval
$[0,\tau_{x_1,x_2})$. 

\begin{Lemma}
\label{l1.2}
$X(t,x_1,x_2)$ has a limit in $L^2(\Omega)$ when $t\to \tau_{x_1,x_2}$. Moreover, 
if $\tau_{x_1,x_2}<T$, we have $X_1(\tau_{x_1,x_2},x_1,x_2)=X_2(\tau_{x_1,x_2},x_1,x_2)$.
\end{Lemma}
{\bf Proof:} 
Let us define 
$$
X_1^\varepsilon(t)= \left\{
\begin{array}{l}
X_1(t,x_1,x_2),\; t\le \tau_{x_1,x_2}^\varepsilon,\\
\\
X(t,\tau_{x_1,x_2}^\varepsilon,X_1(\tau_{x_1,x_2}^\varepsilon,x_1,x_2),\; 
t\ge \tau_{x_1,x_2}^\varepsilon.
\end{array}
\right.
$$
We have denoted by $X(\cdot,s,x)$ the solution of \eqref{e1.1} with the condition $X(s,s,x)=x$
at time $s$ instead of $0$. It is not difficult to check that $X_1$ and $X(\cdot,x_1)$ have 
the same law. Let us write for $\eta_1,\eta_2>0$:
$$
\begin{array}{l}
\E(|X_1(\tau_{x_1,x_2}-\eta_1,x_1,x_2)-X_1(\tau_{x_1,x_2}-\eta_2,x_1,x_2)|^2)\\
\\
=\lim_{\varepsilon\to 0} \E(|X_1(\tau_{x_1,x_2}^\varepsilon-\eta_1)
-X_1(\tau_{x_1,x_2}^\varepsilon-\eta_2)|^2)\\
\\
=\lim_{\varepsilon\to 0} \E(|X_1^\varepsilon(\tau_{x_1,x_2}^\varepsilon-\eta_1)
-X_1^\varepsilon(\tau_{x_1,x_2}^\varepsilon-\eta_2)|^2)\\
\\
\le \lim_{\varepsilon\to 0}\E(\sup_{t\in[0,T]} |X_1^\varepsilon(t-\eta_1,x_1,x_2)
-X_1^\varepsilon(t-\eta_2,x_1,x_2)|^2).
\end{array}
$$
Since, $X_1$ and $X(\cdot,x_1)$ have 
the same law, we can write
$$
\E(\sup_{t\in[0,T]} |X_1^\varepsilon(t-\eta_1)
-X_1^\varepsilon(t-\eta_2)|^2)
= \E(\sup_{t\in[0,T]} |X(t-\eta_1,x_1)
-X(t-\eta_2,x_1)|^2).
$$ 
By \eqref{e1.4bis}, we know that this latter term goes to zero so that we prove that 
$X_1(t)$ has a limit. We treat $X_2(t)$ exactly in the same way.

Finally, if $\tau_{x_1,x_2}<T$, then $|X_1(\tau_{x_1,x_2}^\varepsilon,x_1,x_2)
-X_2(\tau_{x_1,x_2}^\varepsilon,x_1,x_2)|=\varepsilon$ for any $\varepsilon >0$. 
Letting $\varepsilon\to 0$ we deduce the last statement.
\hfill$\Box$\bigskip

We also consider   the following  equation
\begin{equation}
\label{e1.15}
\left\{\begin{array}{lll}
dX_1=(AX_1+b(X_1))dt+\frac{1}{\sqrt{2}}\;\mbox{1{\kern-2.7pt}l}_{[0,\tau_{x_1,x_2}]}(t)dW_1\\
\\
+\frac{1}{\sqrt{2}}\;\mbox{1{\kern-2.7pt}l}_{[0,\tau_{x_1,x_2}]}(t)\left(1-2\frac{(X_1-X_2)\otimes(X_1-X_2)}{|X_1-X_2|^2}\right)dW_2+\frac{1}{\sqrt{2}}
(dW_1+dW_2)(1-\mbox{1{\kern-2.7pt}l}_{[0,\tau_{x_1,x_2}]}(t))\\
\\
dX_2=(AX_2+b(X_2))dt+\frac{1}{\sqrt{2}}\;\mbox{1{\kern-2.7pt}l}_{[0,\tau_{x_1,x_2}]}(t)\left(1-2\frac{(X_1-X_2)\otimes(X_1-X_2)}{|X_1-X_2|^2}\right)dW_1\\
\\
+\frac{1}{\sqrt{2}}\;\mbox{1{\kern-2.7pt}l}_{[0,\tau_{x_1,x_2}]}(t)dW_2+
\frac{1}{\sqrt{2}}(dW_1+dW_2)(1-\mbox{1{\kern-2.7pt}l}_{[0,\tau_{x_1,x_2}]}(t)),\\
\\
X_1(0)=x_1,\quad X_2(0)=x_2.
\end{array}\right. 
\end{equation}
It is clear that for $t\le \tau_{x_1,x_2}$ the solutions of \eqref{e1.5} and \eqref{e1.15} do coincide, whereas for
$t\ge \tau_{x_1,x_2}$ \eqref{e1.15} reduce to
\begin{equation}
\label{e1.16}
\left\{\begin{array}{lll}
dX_1=(AX_1+b(X_1))dt+ \frac{1}{\sqrt{2}}(dW_1(t)+ dW_2(t)),\quad t\ge \tau_{x_1,x_2},\\
\\
dX_2=(AX_2+b(X_2))dt+ \frac{1}{\sqrt{2}}(dW_1(t)+ dW_2(t)),\quad t\ge \tau_{x_1,x_2}.
\end{array}\right. 
\end{equation}
Using Lemma \ref{l1.2}, we easily prove that \eqref{e1.15} has a unique solution.
Moreover, since $\frac{1}{\sqrt{2}}(W_1(t)+W_2(t))$ is 
a cylindrical Wiener process, it follows that $X_1$ and $X_2$ have the same law as 
$X(\cdot,x_1)$ and $X(\cdot,x_2)$. In other words, $(X_1,X_2)$ is a coupling 
of the laws of $X(\cdot,x_1)$ and $X(\cdot,x_2)$.

We are interested in the first time $\tau_{x_1,x_2}$ when $X_1(t,x_1,x_2))$ and $X_2(t,x_1,x_2))$ meet.
 That is $\tau_{x_1,x_2}$ is  the
stopping time
\begin{equation}
\label{e1.9}
\tau_{x_1,x_2}=\inf\{t>0:\;X_1(t,x_1,x_2)=X_2(t,x_1,x_2)\}.
\end{equation}
Our goal is first to show that
\begin{equation}
\label{e1.10}
\E(\tau_{x_1,x_2})<+\infty.
\end{equation}
To prove \eqref{e1.10}  we look, following \cite{CL},  for  a Lyapunov function $f$
such that \eqref{e1.14} holds. This is motivated by next Proposition.
\begin{Proposition}
\label{p2.4}
Assume that there exists a ${\mathcal C}^2$ function such that \eqref{e1.14} holds.
Let   $x_1,x_2\in H$ with $x_1\neq x_2$. Then we have
\begin{equation}
\label{e2.8}
\P(\tau_{x_1,x_2}=+\infty)=0.
\end{equation}
Moreover,
\begin{equation}
\label{e2.9}
\E(\tau_{x_1,x_2})\le   f(|x_1-x_2|).
\end{equation}
\end{Proposition}
  {\bf Proof}. We shall  write for simplicity,
$$
X_1(t)=X_1(t,x_1,x_2),\quad X_1(t)=X_1(t,x_1,x_2).
$$
Also, we can assume without loss of generality that $\kappa=1$.
Then we introduce the following stopping times:
$$
S_N=\inf\left\{t\ge 0: |X_1(t)-X_2(t)|>N\right\},\quad N\in \N,
$$
and
$$
\tau_{n,N}=\tau_{x_1,x_2}^{1/n}\wedge S_N.
$$
By the It\^o formula{\footnote{The application of the It\^o formula can be justified rigorously   
thanks to a regularization argument. This can be done easily in the applications considered 
hereafter.}, Lemma \ref{l1.1} and \eqref{e1.14}, we have
\begin{equation}
\label{e2.11}
\begin{array}{l}
f(|X_1(t\wedge \tau_{n,N})-X_2(t\wedge \tau_{n,N})|) 
\le f(|x_1-x_2|)\\\\
\ds+\int_0^{t\wedge \tau_{n,N}} \frac{f'(|X_1(s)-X_2(s)|)}{|X_1-X_2|}
\tfrac1{\sqrt{2}} (X_1-X_2,d(W_1-W_2))\\\\
- \;(t\wedge \tau_{n,N}).
\end{array}
\end{equation}
It follows  that
$$
\E\big(f(|X_1(t\wedge \tau_{n,N})-X_2(t\wedge \tau_{n,N})|)\big)
\le f(|x-y|)-  \E(t\wedge \tau_{n,N}),
$$
and
$$
\E(t\wedge \tau_{n,N})\le   f(|x-y|) .
$$
Consequently as $t\to \infty$ we find
$$
\E(\tau_{n,N})\le   f(|x-y|)
$$
which yields as $n\to \infty$
$$
\E(\tau_{x_1,x_2}\wedge S_N)\le   f(|x-y|). 
$$
By \eqref{e1.4bis}, we easily prove that $S_N\to \infty$ as $N\to \infty$ so that we get
$$
\E(\tau_{x_1,x_2})\le   f(|x-y|). 
$$
 $\Box$\bigskip

\section{Dissipative systems with white noise}

 We consider the case when there exist $\lambda\ge 0, a>0$ such that
\begin{equation}
\label{e2.1}
(A(x_1-x_2)+ b(x_1)-b(x_2),x_1-x_2   )\le  \lambda|x_1-x_2|^2-a|x_1-x_2|^4,  
\end{equation}
for all $x_1,x_2\in H.$ We also assume that Hypothesis \ref{h1.1} and \eqref{e1.4bis} hold.

A typical equation satisfying such  assumptions is the following
stochastic reaction-diffusion equation on $[0,1]$
$$
\left\{
\begin{array}{l}
dX=(\partial_{\xi\xi} X-\alpha X^3 +\beta	X^2+\gamma X+\delta)dt +dW, \; t>0,\; \xi\in (0,1),\\
\\
X(t,0)=X(t,1)=0,\; t>0,\\
\\
X(0,\xi)=x(\xi), \xi\in (0,1),
\end{array}
\right.
$$
where $\alpha>0$.
In this case, we take $A=D^2_{\xi}$ on the domain $D(A)=H^2(0,1)\cap H^1_0(0,1)$,
$b(x)=-\alpha x^3 +\beta	x^2+\gamma x+\delta$. We could also consider the more
general example where $b$ is a polynomial of degree $2p+1$ with negative leading
coefficient. Note that this equation is gradient, the invariant measure is known explicitly.
However, we shall not use this fact. We could treat as well perturbation of this equation
which are not gradient but satisfy \eqref{e2.1}

Following the above discussion, we look for a  positive function $f$  
such that
\begin{equation}
\label{e2.2}
2\;f''(r)+ f'(r)\left(\lambda r-a\;r^3\right)  =-1.
\end{equation}
  Setting $f'=g$, \eqref{e2.2} becomes
\begin{equation}
\label{e2.3}
g'(r)- \frac12\;g(r) (ar^3-\lambda r)  =-\frac{1}{2}.
\end{equation}
whose general solution is given by
$$
g(r)=e^{\frac{1}{8}\;(ar^4-2\lambda r^2)}g(0)-\frac{1}{2}\int_0^re^{\frac{1}{8}\;(ar^4-as^4-2\lambda r^2+
2\lambda s^2)}ds.
$$
Finally, we have
\begin{equation}
\label{e2.4}
\begin{array}{lll}
f(r)&=&\ds{f(0)+\int_0^re^{\frac{1}{8}\;(as^4-2\lambda s^2)}ds\;f'(0)}\\
\\
&&\ds{-\frac{1}{2}\;\int_0^re^{\frac{1}{8}\;(as^4-2\lambda s^2)}
\left[\int_0^se^{-\frac{1}{8}\;(a\sigma^4-
2\lambda \sigma^2)}d\sigma \right]\;ds.}
\end{array}
\end{equation}
Setting $f(0)=0$ and $f'(0)=\frac{1}{2}\;\int_0^\infty e^{-\frac{1}{8}\;(a\sigma^4-
2\lambda \sigma^2)}d\sigma$ we obtain
\begin{equation}
\label{e2.5}
f(r)=\frac{1}{2}\;\int_0^re^{\frac{1}{8}\;(as^4-2\lambda s^2)}
\left[\int_s^{+\infty}e^{-\frac{1}{8}\;(a\sigma^4-
2\lambda \sigma^2)}d\sigma \right]\;ds 
\end{equation}
and
\begin{equation}
\label{e2.6}
f'(r)=\tfrac{1}{2}\; e^{\frac{1}{8}\;(ar^4-2\lambda r^2)}
 \int_r^{+\infty}e^{-\frac{1}{8}\;(a\sigma^4-
2\lambda \sigma^2)}d\sigma . 
\end{equation}
We need some properties on $f$.  
\begin{Lemma}
\label{l2.1}
We have $(a r^3-\lambda r) f'(r)< 1$ for any $r\ge 0$.
\end{Lemma}
\medskip
{\sf Proof:}    The function $r\mapsto a r^3-\lambda r$
is   increasing and positive if $r>\delta:=\sqrt{\frac{\lambda}{a}}\;$ so that in this case
$$
(a r^3-\lambda r) f'(r)< \tfrac12\;e^{\frac{1}{8}\;(ar^4-2\lambda r^2)}
\int_r^{+\infty} (a \sigma^3-\lambda \sigma) e^{-\frac{1}{8}\;(a\sigma^4-
2\lambda \sigma^2)}d\sigma=1.
$$
Since for $0\le r\le \delta$ we have $(a r^3-\lambda r) f'(r)\le 0$, the conclusion follows. $\Box$\bigskip

\begin{Corollary}
\label{c2.2}
$f' $ is a decreasing positive function.
\end{Corollary}
\medskip
{\sf Proof:} Since $f$ satisfies \eqref{e2.2} , we deduce from Lemma \ref{l2.1} that 
$f''<0$ so that $f'$ decreases. Moreover,  $f'$ is positive by \eqref{e2.6}. $\Box$\bigskip
 
\begin{Lemma}
\label{l2.3}
There exists $\Lambda$ depending only on $a, \lambda$
such that for any $r>0$
\begin{itemize}
\item[i)] $f'(r)\le \Lambda$
\item[ii)] $f(r)\le  \Lambda$
\end{itemize}
\end{Lemma}
\medskip
{\sf Proof:} (i) follows obviously from Corollary \ref{c2.2}. Let us show (ii). Fix
$r>r_0:=\sqrt{\frac{\lambda}{a}}.$
Since $f$ is increasing,   we have
$$
f(r)\le  f(r_0),\quad r\le r_0.
$$
If $r>r_0$ we have by Lemma \ref{l2.1},
$$
\begin{array}{l}
\ds f(r)=   f(r_0)+
\int_{r_0}^r
f'(s)\;ds \\\\
\ds\hphantom{f(r)}\le  f(r_0)+
\int_{r_0}^\infty
\tfrac{ds}{a s^3-\lambda s}=  f(r_0)+ \tfrac1{2\lambda} \ln(1+\tfrac{\lambda}{ar_0^2-\lambda}) .
\end{array}
$$
Therefore  $f(\infty)<\infty$ and ii) follows provided
$
\Lambda\ge \max\{f(\infty),f(r_0)\}.
$
\hfill$\Box$\bigskip

The following results strengthen Proposition \ref{p2.4}.
\begin{Proposition}
\label{p2.5}
Let   $x_1,x_2\in H$ with $x_1\neq x_2$, then there exists $\Lambda>0$ such that
$$
\E\left(e^{{}^{\frac1{2\Lambda^2}\;\scriptstyle\tau_{x_1,x_2}}}\right)\le
e^{\frac1\Lambda}
$$
\end{Proposition}
{\sf Proof:} We use the same notations as in the proof of 
Proposition \ref{p2.4}. By \eqref{e2.11} we have
$$
  t\wedge \tau_{n,N}\le f(|x_1-x_2|)+\int_0^{t\wedge \tau_{n,N}} \frac{f'(|X_1-X_2|)}{|X_1-X_2|}
\tfrac1{\sqrt{2}} (X_1-X_2,d(W_1-W_2)),
$$
so that
\begin{equation}
\label{e2.17}
\E\left( e^{\alpha (t\wedge \tau_{n,N})}\right)\le
e^{ \alpha f(|x_1-x_2|)}\E\left( e^{{}^{ \alpha 
\scriptstyle\int_0^{t\wedge \tau_{n,N}}
\frac{f'(|X_1-X_2|)}{|X_1-X_2|}
\frac1{\sqrt{2}} (X_1-X_2,d(W_1-W_2))}}
\right).
\end{equation}
On the other  hand, since
$$
\E\left( e^{{}^{ \alpha 
\scriptstyle\int_0^{t\wedge \tau_{n,N}}
\frac{f'(|X_1-X_2|)}{|X_1-X_2|}
\frac1{\sqrt{2}} (X_1-X_2,d(W_1-W_2))}}
\right)\le\left(\E\left( e^{{}^{ 2\alpha^2 
\scriptstyle\int_0^{t\wedge \tau_{n,N}}
|f'(|X_1-X_2|)|^2 ds}}
\right)\right)^{1/2},
$$
we deduce from Lemma \ref{l2.3} that
$$
\E\left( e^{{}^{ \alpha 
\scriptstyle\int_0^{t\wedge \tau_{n,N}}
\frac{f'(|X_1-X_2|)}{|X_1-X_2|}
\frac1{\sqrt{2}} (X_1-X_2,d(W_1-W_2))}}
\right)\le\left(\E\left( e^{2\alpha^2 \Lambda^2
(t\wedge \tau_{n,N})}
\right)\right)^{1/2}
$$
Substituting in \eqref{e2.17} yields
 $$
\E\left( e^{\alpha (t\wedge \tau_{n,N})}\right)\le
e^{ \alpha f(|x_1-x_2|)}
\left(\E\left(e^{2\alpha^2 \Lambda^2
(t\wedge \tau_{n,N})}\right)\right)^{1/2}
$$
Choosing $\alpha=\frac1{2\Lambda^2}$, we deduce
$$
\E\left( e^{\frac1{2\Lambda^2} (t\wedge \tau_{n,N})}\right)\le
e^{\frac{1}{2 \Lambda^2}f(|x_1-x_2|)}\le
e^{1/\Lambda}
$$
 Letting $n\to \infty$, $N\to\infty$ and arguing as in the proof of Proposition \ref{p2.4}  we find
the conclusion.\\
\null\hfill$\Box$

\begin{Corollary}
\label{c2.6}
We have
$$
|P_t\varphi(x_1)-P_t\varphi(x_2)|\le
2\;\|\varphi\|_0\;e^{{}^\frac1\Lambda}\;
e^{{}^{-\frac1{2\Lambda^2}t}}
$$
\end{Corollary}
\medskip
{\sf Proof:} Let $x_1,x_2\in H.$ Since $(X_1,X_2)$ is a coupling of $(X(\cdot,x_1),X(\cdot,x_2))$,
we have
$$
|P_t\varphi(x_1)-P_t\varphi(x_2)|=|\E(\varphi(X_1(t))-\varphi(X_2(t)))|
\le \|\varphi\|_0\;\P(\tau_{x_1,x_2}\ge t).
$$
Now the conclusion follows Proposition \ref{p2.5}. $\Box$\bigskip

We end this section by proving that the spectral gap property holds. We thus recover a known 
result (see for instance \cite{DPDG}) with a totally different method. 
\begin{Proposition}
\label{p2.7}
Let $\nu$ be an invariant measure then, for any $p>1$, there exist
$c_p$, $\alpha_p$ such that
$$
|P_t\varphi-\bar\varphi|_{L^p(H,\nu)}\le
c_p e^{-\alpha_p t}\;|\varphi|_{L^p(H,\nu)}
$$
\end{Proposition}
\medskip
{\sf Proof:} By Corollary \ref{c2.6}, we have the result for $p=\infty$.
Using that $P_t$ is a contraction semigroup on $L^1(H,\nu)$
and an interpolation argument, we obtain the result. $\Box$\bigskip

 \section{Burgers equation}
We take here $H=L^2(0,1)$ and denote by $\|\cdot\|$ the norm of the Sobolev space $H^1_0(0,1)$.
We consider the equation
\begin{equation}
\label{e3.1}
\left\{
\begin{array}{l}
dX=(AX+b(X))dt+dW,\\\\
X(0)=x,
\end{array}\right.
\end{equation}
where 
$$
Ax=D^2_\xi x,\quad x\in D(A)=H^2(0,1)\cap H^1_0(0,1)
$$  
  and
$$
b(x)=D_\xi(x^2),\quad x\in H^1_0(0,1).
$$
It well known that problem \eqref{e3.1} has a unique solution for any $x\in L^2(0,1)$ which we denote by
$X(t,x)$, see \cite{DPDT}. It defines a transition semigroup $(P_t)_{t\ge 0}$. It is also known 
that it has a unique invariant measure and is ergodic (see \cite{DPZ2}). The following 
result states the exponential convergence to equilibrium.
\begin{Theorem}
\label{main-burgers}
There exist  constants $C, \gamma >0$ such that for any $x_1,x_2\in L^2(0,1)$, $\varphi\in C_b(L^2(0,1))$,
$$
|P_t\varphi(x_1)-P_t\varphi(x_2)|\le c e^{-\gamma t}\|\varphi\|_0(1+|x_1|^4+|x_2|^4)
$$
\end{Theorem} 
To prove this result,
we want to construct a coupling for equation \eqref{e3.1}. It does not seem possible
to apply directly the method of section 2. We shall first consider a cut off
equation
\begin{equation}
\label{e3.2}
\left\{
\begin{array}{l}
dX=(AX+D_\xi F_R(X))dt+dW,\\\\
X(0)=x,
\end{array}\right.
\end{equation}
where $F_R:L^4(0,1)\to L^2(0,1)$ is defined by
$$
F_R(x)=\left\{\begin{array}{l}
x^2 ,\quad  \mbox{\rm if}\;|x|_{L^4}\le R,\\
\\
\frac{R^2x^2}{|x|^2_{L^4}},\quad  \mbox{\rm if}\;|x|_{L^4}\ge R.
\end{array}\right.
$$
We have 
\begin{equation}
\label{e3.3}
|F_R(x)-F_R(y)| \le 2R\; |x-y|_{L^4},\quad x,y\in L^4(0,1)
\end{equation}
and 
\begin{equation}
\label{e3.4}
|F_R(x)|\le R^2,\quad x\in L^4(0,1).
\end{equation}
We have denoted by $|\cdot|_{L^4}$ the norm in $L^4(0,1)$. The norm in $L^2(0,1)$ is
still denoted by $|\cdot|$.
It is not difficult to check that Hypothesis \ref{h1.1} and \eqref{e1.4bis} hold so that the results
of section 2 can be applied.
We then need a priori estimates on the solutions of \eqref{e3.1} so that we can control
when the coupling for the cut-off equation can be used. These are given in section 4.2.
Then, we construct a coupling for the Burgers equation which enables us to prove 
the result in section 4.4.

\subsection{Coupling for the cut--off equation}
Here $R>0$ is fixed. We denote by the same symbol $c_R$ various constants
depending only on $R$.

\begin{Lemma}
\label{l3.1}
There exists $\alpha>0,\beta>0, c_R>0$ such that
\begin{equation}
\label{e3.5}
\frac{2\alpha}{1-\beta}\in [1,2), 
\end{equation}
and 
\begin{equation}
\label{e3.6}
|F_R(x)-F_R(y)|\le c_R |x-y|^\alpha \|x-y\|^\beta,\quad x,y\in H^1_0(0,1).
\end{equation}
\end{Lemma}
{\bf Proof}. First notice that by \eqref{e3.3} and \eqref{e3.4} it follows that
$$
|F_R(x)-F_R(y)|\le  (2R)^{2-\gamma}\; |x-y|_{L^4}^\gamma,\quad x,y\in L^4(0,1),
$$
for any $\gamma\in[0,1]$. Moreover, by the Sobolev embedding theorem
we have $H^{1/4}(0,1)\subset L^4(0,1)$ and using
 a well known interpolatory inequality we find that
\begin{equation}
\label{e3.6bis}
|x-y|_{L^4}\le c|x-y|_{H^{1/4}}\le c|x-y|^{3/4}\;\|x-y\|^{1/4},\quad x,y\in H^1_0(0,1).
\end{equation}
Consequently
$$
|F_R(x)-F_R(y)|\le  c(2R)^{2-\gamma}\;|x-y|^{3\gamma/4} \|x-y\|^{\gamma/4}
$$
Now setting $\alpha=3\gamma/4,\;\beta=\gamma/4$, the conclusion follows choosing 
 $\gamma\in [\frac47,1]$.
$\Box$\bigskip

We now construct the coupling for equation \eqref{e3.2}.
For any $x,y\in D(A)$ we have, taking into account Lemma \ref{l3.1},
$$
\begin{array}{lll}
(A(x-y)+D_\xi F_R(x)-D_\xi F_R(y),x-y)&\le& -\|x-y\|^2+|F_R(x)-F_R(y)|\;\|x-y\|\\\\
&\le&
-\|x-y\|^2+c_R\;|x-y|^\alpha\;\|x-y\|^{1+\beta}.
\end{array}
$$
Using the elementary inequality
$$
uv\le \tfrac{1+\beta}{2}\;(\epsilon u)^{\frac{2}{1+\beta}}
+\tfrac{1-\beta}{2}\;(\epsilon^{-1} v)^{\frac{2}{1-\beta}},\quad u,v,\epsilon>0,
$$
and choosing suitably $\epsilon$ we find
$$
\begin{array}{lll}
(A(x-y)+D_\xi F_R(x)-D_\xi F_R(y),x-y)&\le&
-\tfrac12 \|x-y\|^2+c_R\;|x-y|^{\frac{2\alpha}{1-\beta}}\\
\\
&\le& -\tfrac12\|x-y\|^2+\tfrac{\pi^2}{4}|x-y|^2+c_R\;|x-y|,
\end{array}
$$
since $\frac{2\alpha}{1-\beta}\in [1,2)$. By the Poincar\'e inequality we conclude that
$$
(A(x-y)+D_\xi F_R(x)-D_\xi F_R(y),x-y)\le
-\tfrac{\pi^2}{4}|x-y|^2+c_R\;|x-y|.
$$
Consequently \eqref{e1.14} (with $\kappa=1$) reduces to
$$
2f_R''(r)-f_R'(r)(\tfrac{\pi^2}4 r-c_R)=-1,
$$
$$
f_R''(r)-f_R'(r)(2\,a\, r-c_R)=-\tfrac12,\qquad a=\tfrac{\pi^2}{16}.
$$
Then
$$
f_R'(r)=e^{a r^2-c_R r}f_R'(0)-\tfrac12 \int_0^r
e^{a (r^2-s^2)-c_R (r-s)}ds
$$
and
$$
f_R(r)=f_R(0)+f_R'(0)\int_0^re^{a s^2-c_R s}ds-\tfrac12 \int_0^r ds\int_0^s
e^{a (s^2-u^2)-c_R (s-u)}du.
$$
  Setting
$$
f_R(0)=0,\quad f_R'(0)=\tfrac12\;e^{a r^2-c_R r}\int_0^{+\infty} 
e^{-a  u^2+c_R u}du,
$$
we obtain the solution
\begin{equation}
\label{enuova}
f_R(r):=\tfrac12\int_0^r e^{as^2-c_Rs}\left(
\int_s^\infty e^{-a\xi^2+c_R\xi}d\xi\right)ds.
\end{equation}

We denote by $X_R(t,x)$ the solution of the cut-off equation \eqref{e3.2}.
The corresponding coupling constructed above is denoted by
$\big(X_{1,R}(t;x_1,x_2),X_{2,R}(t;x_1,x_2)\big)$. Then, setting
$$
\tau_R^{x_1,x_2}=\inf\{t>0 : X_{1,R}(t;x_1,x_2)=X_{2,R}(t;x_1,x_2)\}.
$$
By Proposition \ref{p2.4}, we have
\begin{equation}
\label{e3.8}
\E(\tau_R^{x_1,x_2})\le    f_R(|x_1-x_2|)
\end{equation}
\begin{Remark}
Using similar arguments as in section 3, we can derive bounds on $f_R$ and $f_R'$ and
prove following result for the transition semigroup associated the the cut-off equation. 
For all $\varphi\in C_b(H)$ we have
$$
|P^R_t\varphi(x_1)-P^R_t\varphi(x_2)| \le 2c\;\|\varphi\|_0\;
  (1+|x_1-x_2|)^{1/\Lambda}\;
e^{{}^{-\frac1{2\Lambda^2}t}},\quad x_1,x_2\in H.
$$
\end{Remark}

\subsection{A priori estimates}

Next result is similar to Proposition \ref{p2.4} in \cite{dpd04}.
\begin{Lemma}
\label{l3.7}
Let  $\alpha\ge 0$ and
$$
z^\alpha(t)=\int_{-\infty}^t e^{(A-\alpha)(t-s)} dW(s).
$$
Then, for any $p\in\N$, $\varepsilon>0$, $\delta>0$, there exists a random variable $K(\varepsilon,\delta,p)$ such that
$$
|z^\alpha(t)|_{L^p} \le K(\varepsilon,\delta,p) \;\alpha^{-\frac14+\varepsilon}(1+|t|^\delta)
$$
Moreover, all the moments of $K(\varepsilon,\delta,p)$ are finite.
\end{Lemma}
\medskip
{\sf Proof:} Proceeding as in the proof of Proposition 2.1 in \cite{dpd04}, we have
$$
z^\alpha(t)=\int_{-\infty}^t \Big[(t-\sigma)^{\beta-1}-\alpha
\int_\sigma^t (\tau-\sigma)^{\beta-1} e^{-\alpha(t-\tau)}d\tau\Big]
e^{A(t-\sigma)}Y(\sigma)d\sigma
$$
where
$$
Y(\sigma)=\frac{\sin \pi\beta}{\pi}\int_{-\infty}^\sigma (\sigma-s)^{-\beta}
e^{A(\sigma-s)}dW(s)
$$
and $\beta\in (0,1/4)$. It is proved in \cite{dpd04} that for $\gamma\in [0,1]$
$$
\begin{array}{l}
\ds\Big|(t-\sigma)^{\beta-1}-\alpha
\int_\sigma^t (\tau-\sigma)^{\beta-1} e^{-\alpha(t-\tau)}d\tau\Big| \\\\
\hspace{3cm}\le
c(\beta,\gamma)\;\alpha^{-\gamma}
(t-\sigma)^{\beta-1-\gamma}+
(t-\sigma)^{\beta-1} e^{-\alpha(t-\sigma)}.
\end{array}
$$
We deduce, by Poincar\'e inequality,
$$
|z^\alpha(t)|_{L^p} \le c \int_{-\infty}^t \!\!
\Big(\alpha^{-\gamma}
(t-\sigma)^{\beta-1-\gamma}+
(t-\sigma)^{\beta-1} e^{-\alpha(t-\sigma)}\Big)
e^{-\lambda_p (t-\sigma)} |Y(\sigma)|_{L^p}d\sigma,
$$
and, for $r>1$, $m\in\N$, by H\"older inequality we obtain if $\beta>\frac1{2m}$ and $\beta>\gamma+\frac1{2m}$
$$
\begin{array}{l}
|z^\alpha(t)|_{L^p} \le c \Big[ \alpha^{-\gamma}\Big(\ds\int_{-\infty}^t
(t-\sigma)^{\frac{2m}{2m-1}(\beta-1-\gamma)}(1+|\sigma|^r)^{\frac1{2m-1}}
e^{-\lambda_p(t-\sigma)}d\sigma\Big)^{\frac{2m-1}{2m}}\\\\
\hspace{2.5cm}+\Big(\ds\int_{-\infty}^t
(t-\sigma)^{(\beta-1)\frac{2m}{2m-1}}e^{-\frac{2m\alpha}{2m-1}(t-\sigma)}
(1+|\sigma|^r)^{\frac1{2m-1}}
d\sigma\Big)^{\frac{2m-1}{2m}}  \Big]\\\\
\hspace{6.5cm}\ds\times\Big(\int_{-\infty}^t (1+|\sigma|^r)^{-1} |Y(\sigma)|_{L^p}^{2m}\,d\sigma\Big)^{\frac1{2m}}
\\\\
\hph{|Z^\alpha(t)|_{L^p} }\le c\Big(\alpha^{-\gamma}+\alpha^{-\beta+\frac1{2m}}\Big)(1+|t|^{\frac{r}{2m}})
\Big(\ds \int_{-\infty}^t (1+|\sigma|^r)^{-1} |Y(\sigma)|_{L^p}^{2m}\,d\sigma\Big)^{\frac1{2m}}
\end{array}
$$
and the first statement follows if $\gamma$, $\beta$, $m$ are chosen so that
$$
\tfrac1{2m}< \min(\tfrac{\delta}{r},\tfrac{\varepsilon}2),\qquad
\tfrac14-\tfrac{\varepsilon}2<\gamma+\tfrac{\varepsilon}2<\tfrac14.
$$
Indeed, proceeding as in \cite{dpd04}, we easily prove that
$$
K(\varepsilon,\delta,p)=\Big(\int_{-\infty}^\infty
(1+|\sigma|^r)^{-1} |Y(\sigma)|_{L^p}^{2m}\,d\sigma\Big)^{\frac1{2m}}
$$
has all moments finite.  $\Box$

\bigskip
\begin{Proposition}
\label{p3.8}
Let $x\in L^2(0,1)$ and let $X(t,x)$ be the solution of \eqref{e3.2}.
\begin{itemize}

\item[i)] For any $\delta >0$, there exists a constant $K_1(\delta)$
such that for any $x\in L^4(0,1)$, $t\ge 0$
$$
\E\left(\sup_{s\in [0,t]} |X(s,x)|_{L^4}^4\right)\le 4 |x|_{L^4}^4+K_1(\delta)(1+t^\delta)
$$
\item[ii)] There exists a constant $K_2\ge 0$ such that for any $x\in L^4(0,1)$
and $t\ge 0$
$$
\E( |X(t,x)|_{L^4}^4)\le (e^{-\pi^2 t/16}|x|_{L^4}+K_2)^4.
$$
\item[iii)] There exists a constant $K_3$ such that for any $x\in L^2(0,1)$
and $t\in[1,2]$
$$
\E( |X(t,x)|_{L^4}^4)\le K_3 (1+|x|_{L^2}^4).
$$
\end{itemize}
\end{Proposition}
\medskip
{\sf Proof:} In the proof we shall  denote by $c$ several different constant. Let us prove (i).
Fix $t>0, x\in L^4$ and $\delta>0$ and set
$Y(s)=X(s,x)-z^\alpha(s)$,
 $Y(s)$   satisfies
$$
\left\{\begin{array}{l}
\ds\frac{dY(s)}{ds}=AY(s)+D_\xi(Y(s)+z^\alpha(s))^2+\alpha z^\alpha(s),\\\\
Y(0)=x-z^\alpha(0)
\end{array}\right.
$$
By  similar computations as in  \cite[Proposition 2.2]{dpd04}, we have   that
$$
\tfrac14\; \frac{d\ }{ds}|Y(s)|_{L^4}^4+\tfrac32 \int_0^1|Y(s)|^2\,|D_\xi Y(s)|^2\,d\xi
\le c |z^\alpha(s)|_{L^4}^{\frac83}\,|Y(s)|_{L^4}^4+
c |z^\alpha(s)|_{L^4}^4+\alpha^4|z^\alpha(s)|_{L^4}^4,
$$
and, by the Poincar\'e inequality,
\begin{equation}
\label{e1}
\frac{d\ }{ds} |Y(s)|_{L^4}^4+\tfrac{3\pi^2}8 |Y(s)|_{L^4}^4\le
c |z^\alpha(s)|_{L^4}^{\frac83} |Y(s)|_{L^4}^4+(c+\alpha^4(s))|z^\alpha(s)|_{L^4}^4.
\end{equation}
We now choose $\alpha$ so large that
\begin{equation}
\label{e2}
c |z^\alpha(s)|_{L^4}^{\frac83}\le \tfrac{\pi^2}8,\quad s\in [0,t].
\end{equation}
For this we use Lemma \ref{l3.7} with
$$
\varepsilon=\tfrac18,\quad p=4 \quad\mbox{\rm and}\;\delta\;\mbox{\rm replaced by}\;\tfrac\delta 8.
$$ 
We see that there exists $ c>0$ such that \eqref{e2} holds provided
$\alpha= c \;\big( K(\tfrac18,\tfrac{\delta}8,4)(1+t^{\delta/8}\big)^8$.
So, by \eqref{e1} it follows that
$$
\frac{d\ }{ds}|Y(s)|_{L^4}^4+\tfrac{\pi^2}4 |Y(s)|_{L^4}^4\le
c\,\Big(1+K(\tfrac18,\tfrac{\delta}8,4)^8(1+t)^\delta\Big),\quad s\in [0,t].
$$
Consequently, by the Gronwall lemma, we see that
$$
\begin{array}{l}
|Y(s)|_{L^4}^4\le e^{-\frac{\pi^2}4s}|Y(0)|_{L^4}^4+c\ds\int_0^s
e^{-\frac{\pi^2}4(s-\sigma)}(1+K(\tfrac18,\tfrac{\delta}8,4)^8(1+\sigma)^\delta d\sigma\\\\
\hph{|Y(s)|_{L^4}^4}\le e^{-\frac{\pi^2}4s} |Y(0)|_{L^4}^4+c(1+K(\tfrac18,\tfrac{\delta}8,4)^8(1+t^\delta)),
\quad s\in [0,t],
\end{array}
$$
which yields by \eqref{e2}
\begin{equation}
\label{e3}
\begin{array}{l}
|X(s,x)|_{L^4}^4\le 4e^{-\frac{\pi^2}4s}|x|_{L^4}^4+4|z^\alpha(0)|^4_{L^4}+4|z^\alpha(s)|^4_{L^4}
+c\,(1+K(\tfrac18,\tfrac{\delta}8,4)^8(1+t^\delta) )\\\\
\hph{|X(s,x)|_{L^4}^4}\le 
4e^{-\frac{\pi^2}4s}|x|_{L^4}^4+c(1+K(\tfrac18,\tfrac{\delta}8,4)^8(1+t^\delta)),
\quad s\in [0,t].
\end{array}
\end{equation}
Now, i) follows since $K(\tfrac18,\tfrac{\delta}8,4)$ has  finite moments.\bigskip

To prove ii) we denote by $X(t,-t_0;x)$ the solution at time $t$ of the Burgers equation
with initial data $x$ at the time $-t_0$.
 Since $X(t_0,x)$ and $X(0,-t_0;x)$ have the same
law, it suffices to prove
$$
\E(|X(0,-t_0;x)|_{L^4}^4)\le (e^{-\pi^2t_0/16} |x|_{L^4}+K_2)^4,\quad t_0\ge 0.
$$
We set $Y(t)=X(t,-t_0:x)-z^\alpha(t)$. 
Proceeding as above we find
$$
|Y(t)|_{L^4}^4\le e^{-\frac{\pi^2}4(t-s)}|Y(s)|_{L^4}^4+c(1+K(\tfrac18,\tfrac{1}8,4)^8(1+|s|)),
\mbox{ for } 0\ge t\ge s.
$$
Since $(a+b)^{1/4}\le a^{1/4}+b^{1/4}$, we obtain
$$
|X(t,-t_0;x)|_{L^4}\le e^{-\frac{\pi^2}{16}(t-s)}|X(s,-t_0;x)|_{L^4}
+c(1+K(\tfrac18,\tfrac{1}8,4)^8(1+|s|))^{1/4},
\mbox{ for } 0\ge t\ge s.
$$
We choose $-s=-t+1=1,2,\dots,n_0$ with $n_0=[t_0]$ and then  $-s=-t+1=t_0$, we obtain
$$
\begin{array}{ll}
|X(0,-t_0,x)|_{L^4} &\le e^{-n_0\pi^2(t-s)/16} |X(-n_0,-t_0,x)|_{L^4}\\
\\
&+c 
\sum_{l=0}^{n-1}e^{-l\pi^2/16}  (1+K(\tfrac18,\tfrac{1}8,4)^8(1+l))^{1/4}\\
\\
&\le e^{-t_0\pi^2(t-s)/16} |x|_{L^4}+c 
\sum_{l=0}^{n-1}e^{-l\pi^2/16}  (1+K(\tfrac18,\tfrac{1}8,4)^8(1+l))^{1/4}\\
\\
&+ e^{-n_0\pi^2/16}  (1+K(\tfrac18,\tfrac{1}8,4)^8(1+t_0))^{1/4}\\
\\
&\le e^{-t_0\pi^2(t-s)/16} |x|_{L^4}+c (1+K(\tfrac18,\tfrac{1}8,4)^8)^{1/4},
\end{array}
$$
so, ii) follows.\bigskip

To prove iii), we use, as in i), $Y(s)=X(s,x)-z^\alpha(s)$ and have
$$
\tfrac12\;\frac{d\ }{ds}|Y(s)|_{L^2}^2+|D_\xi Y(s)|_{L^2}^2\le
c|z^\alpha(s)|_{L^4}^{\frac83} |Y(s)|_{L^2}^2+c|z^\alpha(s)|_{L^4}^4 +\alpha^2|z^\alpha(s)|_{L^2}^2
$$
so that, choosing $\alpha$ conveniently,
$$
\sup_{[0,2]}|Y(s)|_{L^2}^2+\int_0^2 |D_\xi Y(s)|_{L^2}^2 ds\le
 |x|_{L^2}^2 +c(1+\alpha^2).
$$
For instance, we can take $\alpha=c K(\tfrac18,1,4)^8$.\\
Moreover
$$
\frac{d\ }{ds}\big( s|Y(s)|_{L^4}^4\big)+\tfrac{\pi^2}4 s |Y(s)|_{L^4}^4\le |Y(s)|_{L^4}^4
+c(1+\alpha^4)
$$
so that for $t\in[1,2]$
$$
|Y(t)|_{L^4}^4 \le \int_0^2 |Y(s)|_{L^4}^4 ds+c(1+\alpha^4).
$$
Using the inequality \eqref{e3.6bis}, we have
$$
|Y(t)|_{L^4} \le c\; |Y(t)|_{L^2}^{\frac34} \;|D_\xi Y(t)|_{L^2}^{\frac14}
$$
and we deduce
$$
\begin{array}{l}
|Y(t)|_{L^4}^4\le c\ds\int_0^2|Y(s)|_{L^2}^3|D_\xi Y(s)|_{L^2}ds +c(1+\alpha^4)\\\\
\hph{|Y(t)|_{L^4}^4}\le c(|x|_{L^2}^4 +(1+\alpha^2)^2)+c(1+\alpha^4)
\end{array}
$$
for $t\in[1,2]$ and iii) follows.
 $\Box$
\bigskip

Next lemma is similar to Lemma 2.6 in Kuksin-Shirikyan \cite{ks}.

\begin{Lemma}\label{l3.9}
Let $\rho_0>0$ and $\rho_1>0$, there exist $\alpha(\rho_0,\rho_1)>0$ and $T(\rho_0,\rho_1)>0$
such that for $t\in[T(\rho_0,\rho_1),2T(\rho_0,\rho_1)]$
and $|x_1|_{L^4}\le \rho_0$, $|x_2|_{L^4}\le \rho_0$,
$$
\P\big(|X(t,x_1)|_{L^4} \le \rho_1
\mbox{ and }
|X(t,x_2)|_{L^4} \le \rho_1\big)\ge \alpha(\rho_0,\rho_1).
$$
\end{Lemma}
\medskip
{\sf Proof:} Let $X^0$ be the solution of the deterministic Burgers equation
$$
\left\{
\begin{array}{l}
\ds\frac{dX^0}{dt}(t,x)=AX^0(t,x) + b(X^0(t,x))\\\\
X^0(0,x)=x.
\end{array}
\right.
$$
Since, as easily checked, $(b(X^0),(X^0)^3)=0$, we obtain by standard computations
$$
|X^0(t,x)|_{L^4}^4\le e^{-\pi^2t/4 } |x|_{L^4}^4\le e^{-\pi^2t/4 } \rho_0^4,\quad\mbox{\rm for all}\;|x|_{L^4}\le \rho_0.
$$
We choose
$$
T(\rho_0,\rho_1)=\frac{16}{\pi^2}\; \ln\Big(\frac{2\rho_0}{\rho_1}\Big).
$$
Then, for $t\ge T(\rho_0,\rho_1)$
$$
|X^0(t,x)|_{L^4} \le \frac{\rho_1}2\quad\mbox{\rm for all}\;|x|_{L^4}\le \rho_0.
$$
let
$$
z(t)=\int_0^t e^{(t-s)A}dW(s)
$$
then for any $\eta>0$
$$
\P\Big(\sup_{[0,2T_0(\rho_0,\rho_1)]} |z(t)|_{L^4} \le \eta\Big) >0
$$
and since the solution of the stochastic equation is a continuous function
of $z$ we can find $\eta$ such that
$$
|X(t,x_i)|_{L^4} \le |X^0(t,x_i)| +\frac{\rho_1}{2},\; i=1,2,
$$
for $t\in [0,2 T(\rho_0,\rho_1)]$ provided $|z(t)|_{L^4} \le \eta$ on 
$[0,2 T(\rho_0,\rho_1)]$.
The conclusion follows easily. $\Box$
\bigskip

\subsection{Construction of the coupling}

Let $x_1,x_2\in L^4(0,1)$, we construct
$\big(X_1(\cdot;x_1,x_2),X_2(\cdot;x_1,x_2)\big)$ a coupling of $X(\cdot;x_1)$ and $X(\cdot;x_2)$  as follows. 
Fix   $\rho_0>0$, $\rho_1>0$,
$R>\max\{\rho_0,\rho_1\}$, $T>T_0:=T(\rho_0,\rho_1)$ (defined in Lemma \ref{l3.9}), all  to be chosen later. 

We recall that $\big(X_{1,R}(t;x_1,x_2),X_{2,R}(t;x_1,x_2)\big)$ represents the coupling of $X_R(\cdot;x_1)$ and
$X_R(\cdot;x_2)$ constructed in section 4.1,   
where $X_R(\cdot;x_1)$ and $X_R(\cdot;x_2)$  are the solutions of the cut-off equation
\eqref{e3.2}.

We shall need also the  coupling  of $X_R(\cdot,T_0;x_1)$ and $X_R(\cdot,T_0;x_2)$   when the initial time is any $T_0>0$ instead of 0.
We denote it by
$$
\big(X_{1,R}(t,T_0;x_1,x_2),X_{2,R}(t,T_0;x_1,x_2)\big)
$$
and in this case we shall write
$$
\tau_{R}^{x_1,x_2}=\inf\{t>0 : X_{1,R}(t+T_0,T_0;x_1,x_2)=X_{2,R}(t+T_0,T_0;x_1,x_2)\}.
$$
Notice that  $\tau_{R}^{x_1,x_2}$ does not  depend on $T_0$ thank to the Markov property because the
Burgers equation does not depend explicitely on time.\bigskip

\noindent First we shall   construct the coupling on $[0,T]$,   defining $(X_1(t;x_1,x_2),X_2(t;x_1,x_2))$ as follows.
If 
\begin{equation}
\label{condition}
x_1\not = x_2,\; |x_1|_{L^4}\le \rho_0,\; |x_2|_{L^4}\le \rho_0,\; |X(T_0;x_1)|_{L^4}\le \rho_1\mbox{ and }|X(T_0;x_2)|_{L^4}\le \rho_1,
\end{equation}
we set
$$
X_1(t;x_1,x_2)=X(t,x_1),\;\; X_2(t;x_1,x_2)=X(t,x_2),\quad\mbox{\rm for}\;t\in [0,T_0]
$$
and for $i=1,2$
$$
X_i(t;x_1,x_2)=\left\{\begin{array}{l}
X_{i,R}(t,T_0;X(T_0,x_1),X(T_0,x_2))\;
\mbox{\rm if}\;T_0\le t\le  \min\{\widetilde{\tau}_R,T\},\\\\
X(t,\tilde\tau_R;X_{i,R}(\tilde\tau_R,T_0;X(T_0,x_1),X(T_0,x_2))\;
\mbox{\rm if}\; \min\{\widetilde{\tau}_R,T\}<t\le T,
\end{array}\right.
$$
where
$$
\widetilde{ \tau}_R=\inf\{t\ge T_0: \max\{|X_{i,R}(t,T_0;X(T_0,x_1),X(T_0,x_2))|_{L^4},\;i=1,2\}>R\}.
$$
If \eqref{condition} does not hold, we simply set 
$$
X_1(t;x_1,x_2)=X(t,x_1),\;\; X_2(t;x_1,x_2)=X(t,x_2),\quad\mbox{\rm for}\;t\in [0,T]
$$
So, we have constructed the coupling on $[0,T]$.
 The preceding construction can be obviously generalized considering a time interval $[t_0,t_0+T]$
and  random initial data $(\eta_1,\eta_2)$,  $\mathcal F_{t_0}$-measurable. In this case we
denote the coupling by
$$
(X_1(t,t_0,\eta_1,\eta_2),X_2(t,t_0,\eta_1,\eta_2)).
$$
Now we define the coupling
$(X_1(t;x_1,x_2),X_2(t;x_1,x_2))$ for all time, setting by recurrence
$$
X_i(t,x_1,x_2):=X_i(t,kT,X_1(kT,x_1,x_2),X_2(kT,x_1,x_2)),\quad i=1,2,\;t\in [kT,(k+1)T].
$$
Let us summarize the construction of the coupling on $[0,T]$. We first let the  original
processes $X(\cdot,x_1)$, $X(\cdot,x_2)$ evolve until they are both
in the ball of radius $\rho_0$. Then, we let them evolve and if at time $T_0$ they both enter the
ball of radius $\rho_1$ we use the coupling of the truncated equation as long
as the norm do not exceed $R$ (so that if $\widetilde{\tau}_R\ge T-T_0$ we have
a coupling of the Burgers equation having good properties).
Then, if the coupling is successful, i.e. if $X_1(T;x_1,x_2)=X_2(T;x_1,x_2)$,
we use the original Burgers equation and the solutions remain equal. Otherwise,
we try again in $[T,2T]$ and so on.

\subsection{Exponential convergence to equilibrium for the Burgers equation}

We shall choose  now $\rho_0$, $\rho_1$, $R$ and $T$ (recall that $T_0=T(\rho_0,\rho_1)$ is determined
by Lemma \ref{l3.9}).   
We first assume that
$$
|x_1|_{L^4}\le \rho_0\qquad \mbox{and}\qquad |x_2|_{L^4}\le \rho_0
$$
and set
$$
\begin{array}{l}
\ds A=\left\{\tau_R^{X(T_0,x_1),X(T_0,x_2)}\le T-T_0\right\}\quad(\mbox{\rm the coupling is successfull in }\;[T_0,T]),\\
\\
\ds{B=\left\{\sup_{t\in [T_0,T],i=1,2}|X_{i,R}(t,T_0;X(T_0,x_1),X(T_0,x_2))|_{L^4}
\le R\right\},}\\
\\
\ds{C=\left\{\max_{ i=1,2} |X (T_0,x_i)|_{L^4}\le \rho_1 \right\}.}
\end{array}
$$
Then, we have
$$
\P\big( X_1(T;x_1,x_2)=X_2(T,x_1,x_2)\big)\ge\P(A\cap B\cap C).
$$
We are now going to estimate $\P(A\cap B\cap C)$.
Concerning $C$ we note that by Lemma \ref{l3.9} it follows that
$$
\P\left(\max_{ i=1,2}|X(T_0,x_i)|_{L^4}\le \rho_1\right)\ge \alpha(\rho_0,\rho_1).
$$
Moreover,  
\begin{equation}
\label{e3.9}
\begin{array}{l}
\P(A\cap B\cap C)=\\Ê
\\
\ds{\hspace{-.8cm}\int\limits_{\hspace{1cm}|y_i|_{L^4}\le \rho_1,i=1,2}\hspace{-1.2cm}\P(A\cap B|
X(T_0,x_1)=y_1,X(T_0,x_2)=y_2)\P(X(T_0,x_1)\in dy_1,X(T_0,x_2)\in dy_2).}
\end{array}
\end{equation}
But
\begin{equation}
\label{e3.10}
\begin{array}{l}
\P(A\cap B|X(T_0,x_1)=y_1,X(T_0,x_2)=y_2)\\
\\
=\P(\tau_R^{y_1,y_2} \le T-T_0 \mbox{ and }
\ds\sup_{[T_0,T]} |X_{i,R}(t,T_0;y_1,y_2)|_{L^4}\le R\quad \mbox{for } i=1,2)\\\\
\ge 1-\P(\tau_R^{y_1,y_2}> T-T_0)-\ds\sum_{i=1,2}
\P\left(\sup_{[T_0,T]} |X_{i,R}(t,T_0;y_1,y_2)|>R\right)
\end{array}
\end{equation}
By  the Chebyshev inequality and \eqref{e3.8} it follows that
$$
\P(\tau_R^{y_1,y_2}\ge T-T_0)\le \frac1{T-T_0}\;\E(\tau_R^{y_1,y_2})\le  \frac1{T-T_0}\;f ^R(|y_1-y_2|)
$$
and, since $\mathcal L\big(X_{i,R}(\cdot,T_0;y_1,y_2)\big)=
\mathcal L\big(X_R(\cdot,T_0;y_i)\big)=\mathcal L\big(X_R(\cdot;y_i)\big)$, we have, taking into account Proposition
\ref{p3.8}-(i), that
$$
\begin{array}{l}
\ds\P\left(\sup_{[T_0,T]}|X_{i,R}(t,T_0;y_1,y_2)|_{L^4}> R\right)=
\P\left(\sup_{[0,T-T_0]}|X_R(t,y_i)|_{L^4}> R\right)\\
\\ 
\ds =\P\left(\sup_{[0,T-T_0]}|X(t,y_i)|_{L^4}> R\right) 
\ds\le \frac{4\rho_1^4+K_1(\delta)(1+(T-T_0)^\delta)}{R^4}.
\end{array}
$$
 Consequently, if $|y_i|_{L^4}\le \rho_1$ for $i=1,2$,
$$
\begin{array}{l}
\ds 1-\P(\tau_R^{y_1,y_2}> T-T_0)-\sum_{i=1,2}
\P\left(\sup_{[T_0,T]} |X_{i,R}(t,T_0;y_1,y_2)|>R\right)\\\\
\ds \hspace{3cm}\ge 1-\frac1{T-T_0}\; f_R(|y_1-y_2|)-
2\frac{4\rho_1^4+K_1(\delta)(1+(T-T_0)^\delta)}{R^4}\\\\
\ds \hspace{3cm}\ge 1-\frac1{T-T_0}\; f_R(2\rho_1)-
\frac{8\rho_1^4+2K_1(\delta)(1+(T-T_0)^\delta)}{R^4}
\end{array}
$$
We deduce by \eqref{e3.9} and \eqref{e3.10} that for $|x_i|_{L^4}\le \rho_0$, $i=1,2$,
\begin{equation}
\label{e3.11}
\begin{array}{l}
\P\big(X_1(T,x_1,x_2)=X_2(T,x_1,x_2)\big)\ge \big(1-\frac1{T-T_0}\;f_R(2\rho_1)\\
\\
-\frac{8\rho_1^4+2K_1(\delta)(1+(T-T_0)^\delta)}{R^4}\big) \alpha(\rho_0,\rho_1).
\end{array}
\end{equation}
We choose now $T-T_0=1$, $\rho_1\le 1$ and $R$ such that
\begin{equation}
\label{e3.12}
\frac{8+4K_1(\delta)}{R^4}\le \frac14.
\end{equation}
Then, we take $\rho_1$ such that
\begin{equation}
\label{e3.13}
 f_R(2\rho_1)\le \frac14
\end{equation}
(this is possible since $f_R(0)=0$ and $f_R$ is continuous).\\
It follows
\begin{equation}
\label{e3.14}
\P\big(X_1(T;x_1,x_2)=X_2(T;x_1,x_2)\big)\ge \frac12\;\alpha(\rho_0,\rho_1),\quad\mbox{\rm for all }\;
|x_1|_{L^4}\le \rho_0,\;|x_2|_{L^4}\le \rho_0.
\end{equation}
\bigskip

To treat the case of arbitrary $x_1$, $x_2$,  have to choose $\rho_0$. We proceed as in
Kuksin-Shirikyan \cite{ks} and introduce the following
Kantorovich functional
$$
F_k=\E\big( (1+\nu(|X_1(kT;x_1,x_2)|_{L^4}^4+|X_2(kT;x_1,x_2)|_{L^4}^4))
\mbox{1{\kern-2.7pt}l}_{X_1(kT;x_1,x_2)\neq X_2(kT;x_1,x_2)} \big),
$$
where $k\in \N\cup \{0\}$ and $\nu$ is to be chosen later.
\begin{Proposition}
\label{t3.11}
There exist positive numbers $\rho_0,\nu,\gamma$ such that   
\begin{equation}
\label{t3.15}
\P(X_1(kT;x_1,x_2)\neq X_2(kT;x_1,x_2))\le e^{-\gamma k}
\;(1+\nu (|x_1|_{L^4}^4+|x_2|_{L^4}^4)),\quad x_1,x_2\in L^4(0,1).
\end{equation}
\end{Proposition}
{\bf Proof}. We shall denote by the same symbol $c$ several different constants.
Let us estimate $F_1$ in terms of $F_0=(1+\nu(|x_1|_{L^4}^4+|x_2|_{L^4}^4))\mbox{1{\kern-2.7pt}l}_{x_1\neq x_2}$.
If
$x_1=x_2$ then $X_1(T;x_1,x_2)=X_2(T;x_1,x_2)$ a.s. and so, $F_1=0$.
Let now $x_1\neq x_2$. If $|x_1|_{L^4}>\rho_0$ then $X_1(T;x_1,x_2)=X(T,x_1)$ and $X_2(T;x_1,x_2)=X(T,x_2)$.
Consequently, taking into account Proposition \ref{p3.8}--(ii),
$$
\begin{array}{l}
F_1=\E\left[(1+\nu(|X(T,x_1)|_{L^4}^4+|X(T,x_2)|_{L^4}^4))
\mbox{1{\kern-2.7pt}l}_{X(T,x_1)\neq X(T,x_2)}\right]\\\\
\hspace{1cm}\le \E\big[1+\nu(|X(T,x_1)|_{L^4}^4+|X(T,x_2)|_{L^4}^4)\big]\\\\
\hspace{1cm}\le 1+\nu\big((e^{-\pi^2T/16} |x_1|_{L^4}+K_2)^4+
(e^{-\pi^2T/16} |x_2|_{L^4}+K_2)^4 \big).
\end{array}
$$
Since $T\ge T-T_0=1$, there exists $c$ such that
\begin{equation}
\label{t3.16}
(e^{-\pi^2T/16} a+b)^4\le e^{-\pi^2T/8} a^4+c \;b^4,
\qquad \mbox{for any } a,b\ge 0
\end{equation}
and we deduce that, if  $x_1\neq x_2$ and $|x_1|_{L^4}>\rho_0$ (or if $x_1\neq x_2$ and $|x_2|_{L^4}>\rho_0$),
\begin{equation}
\label{t3.17}
F_1\le 1+\nu\big(e^{-\pi^2T/8}
(|x_1|_{L^4}^4+|x_2|_{L^4}^4)+cK_2^4 \big)
\end{equation}
Let us now consider  the case when $x_1\neq x_2$, $|x_1|_{L^4}\le \rho_0$ and $|x_2|_{L^4}\le \rho_0$.
Taking into account \eqref{e3.14} and Proposition \ref{p3.8}, we have
\begin{equation}
\label{t3.18}
\begin{array}{l}
F_1=\E\big[(1+\nu(|X_1(T;x_1,x_2)|_{L^4}^4+|X_2(T;x_1,x_2)|_{L^4}^4))
\mbox{1{\kern-2.7pt}l}_{X_1(T;x_1,x_2)\neq X_2(T;x_1,x_2)}\big]\\\\
\le \P(X_1(T;x_1,x_2)\neq X_2(T;x_1,x_2))+\nu
\E[|X_1(T;x_1,x_2)|_{L^4}^4+|X_2(T;x_1,x_2)|_{L^4}^4]\\\\
\ds \le 1-\frac12\; \alpha(\rho_0,\rho_1)+\nu \big(e^{-\pi^2T/8}
(|x_1|_{L^4}^4+|x_2|_{L^4}^4)+cK_2^4\big)
\end{array}
\end{equation}
since $\mathcal L(X_i(T;x_1,x_2))=\mathcal L(X(T;x_i))$ and thanks to
Proposition \ref{p3.8}--ii) and \eqref{t3.16}.

To conclude,  we  shall choose $\rho_0$ and $\nu$   such that
\begin{equation}
\label{t3.18bis}
\begin{array}{l}
q_1(\lambda):=
\ds\frac{1+\nu(e^{-\pi^2 T/8}\lambda+
cK_2^4)}{1+\nu\lambda}\le e^{-\gamma},\quad\mbox{\rm for }\lambda> \rho_0^4\\\\
q_2(\lambda):=
\ds\frac{1-\frac12\; \alpha(\rho_0,\rho_1)+\nu(e^{-\pi^2 T/8}\lambda+
cK_2^4)}{1+\nu\lambda}\le e^{-\gamma},\quad\mbox{\rm for }\lambda\le  2\rho_0^4.
\end{array}
\end{equation}
Note first  that $q_1$ is decreasing in $\lambda$ and  
tends to $e^{-\pi^2 T/8}$ as $\lambda \to \infty$. 
We choose $\rho_0$ such that
$$
\rho_0^4= \frac{2(1+ cK_2^4)}{1- e^{-\pi^2T/8}},\quad T=1+T_0(\rho_0,\rho_1).
$$
(It is easy to see that this equation has
  a solution). With this choice of
$\rho_0$, it is also easy to check that $q_1(\rho^4_0)<1$.\\
Choosing now
$$
\nu = \frac{\alpha(\rho_0,\rho_1)}{4cK_2^4}
$$
it is easy to check that
$$
q_2(\lambda)\le \max\{e^{-\pi^2T/8},1-\tfrac14\; \alpha(\rho_0,\rho_1)\}.
$$
Hence,  choosing $\gamma$ such that
$$
e^{-\gamma} \ge \max\{1-\tfrac14\; \alpha(\rho_0,\rho_1),q_1(\rho_0)\},
$$
\eqref{t3.18bis} is fulfilled.

Now we continue the estimate of $F_1$. For any $\lambda\ge \rho_0^4$ we have from the first inequality in \eqref{t3.18bis}
  that
$$
1+\nu( e^{-\pi^2T/8} \lambda+c K_2^4)\le e^{-\gamma}(1+\nu \lambda).
$$
Then  from \eqref{t3.17} we deduce that if $x_1\neq x_2$ and $|x_1|_{L^4}>\rho_0$ or 
$|x_2|_{L^4}>\rho_0$, we have
$$
F_1\le e^{-\gamma}(1+\nu (|x_1|_{L^4}^4+|x_2|_{L^4}^4 ))=e^{-\gamma} F_0.
$$
Moreover
$$
1-\frac12\; \alpha(\rho_0,\rho_1)+\nu \big(e^{-\pi^2T/8} \lambda+cK_2^4\big)
\le e^{-\gamma}(1+\nu \lambda)
$$
for any $\lambda \le 2 \rho_0^4$. Then, by \eqref{t3.18}, we obtain 
for $x_1\neq x_2$, $|x_1|_{L^4}\le \rho_0$ and $|x_2|_{L^4}\le \rho_0$
$$
F_1\le e^{-\gamma}(1+\nu (|x_1|_{L^4}^4+|x_2|_{L^4}^4 ))=e^{-\gamma} F_0.
$$
Therefore, we have in any case
$$
F_1\le e^{-\gamma}\;F_0,\quad x_1,x_2\in L^4(0,1).
$$
It is not difficult to check that $(X_1(kT;x_1,x_2),X_2(kT;x_1,x_2))_{k\in\N}$
is a Markov chain so that we obtain for any $k\in\N$
$$
F_{k+1}\le e^{-\gamma}\;F_k,\quad x_1,x_2\in L^4(0,1).
$$
and, so
$$
F_{k}\le e^{-k\gamma}\;F_0,\quad x_1,x_2\in L^4(0,1).
$$
In particular
$$
\P(X_1(kT;x_1,x_2)\neq X_2(kT;x_1,x_2))\le
e^{-k\gamma}\;(1+\nu (|x_1|_{L^4}^4+|x_2|_{L^4}^4)),\quad x_1,x_2\in L^4(0,1).
$$
$\Box$\bigskip

By Proposition \ref{t3.11} the exponential convergence to equilibrium follows for
$x_1,x_2\in L^4(0,1)$.
If $x_1,x_2\in L^2(0,1)$, we write
$$
\begin{array}{l}
\P(X_1(kT;X(1,x_1),X(1,x_2))\neq X_2(kT;X(1,x_1),X(1,x_2))\\\\
\hspace{2cm}\le 
 e^{-(k-1)\gamma}\;\big(1+\nu \E(|X(1,x_1)|_{L^4}^4+|X(1,x_2)|_{L^4}^4)\big) 
\end{array}
$$
and use Proposition \ref{p3.8}--(iii) to conclude the proof of Theorem \ref{main-burgers}. $\Box$


\begin{thebibliography}{99}

\bibitem{CL} M. F. Chen, S. F. Li,  {\it Coupling methods for multidimensional diffusion processes},
The Annals of Probability, {\bf 17}, no. 1, 151-177,  1989.

\bibitem{dpd04} G. Da Prato, A. Debussche, {\it $m$-dissipativity
of Kolmogorov operators corresponding to Burgers equations with
space-time white noise}, Preprint SNS. Pisa, 2004.

\bibitem{DPDG} G. Da Prato, A. Debussche and B.
Goldys, {\it Invariant measures of non symmetric dissipative stochastic systems,} Probab. Theory Relat. Fields, 
{\bf 123}, 3, 355-380.



\bibitem{DPDT} G. Da Prato, A. Debussche and R. Temam {\it Stochastic Burgers equation,}  NoDEA, 389--402, 1994. 

\bibitem{DPZ2} G. Da Prato and J. Zabczyk, {\it Ergodicity for Infinite Dimensional Systems,} London
Mathematical Society Lecture Notes, n.229, Cambridge University Press, 1996.


\bibitem{Hairer1} M. Hairer, {\it Exponential mixing properties of stochastic PDEs through asymptotic coupling,}
Probab. Theory Related Fields {\bf 124}, 345--380,  2002.


\bibitem{KS1}  S. Kuksin and A. Shirikyan,
{\it Ergodicity for the randomly forced 2D Navier-Stokes equations}, Math. Phys. Anal. Geom. {\bf 4}, 2001.


\bibitem{ks} S. Kuksin S., A. Shirikyan, {\it Coupling approach to white-forced
nonlinear PDEs},  J. Math. Pures Appl. {\bf 81}, 567-602, 2002


\bibitem{LR} T. Lindvall, L. C. G. Rogers, {\it Coupling of mulidimensional diffusions
by reflection}, Ann. Prob., {\bf 14}, 860-872, 1986.


\bibitem{mattingly1} J.C. Mattingly,
{\it Ergodicty of 2D Navier-Stokes equations with random forcing and large viscosity}, 
Commun. Math. Phys. {\bf 206}, 273--288, 1999.



\bibitem{mattingly2} J.C. Mattingly,
{\it Exponential convergence for the stochastically forced Navier-Stokes equations and other partially
dissipative dynamics},  Preprint, 2001.


\bibitem{mueller}  C. Mueller, {\it Coupling and invariant measures for the 
 heat equation with noise,} Ann. Prob. {\bf 21}, 2189--2199, 1993.


\bibitem{odasso} C. Odasso, {\it Ergodicity for  the stochastic complex Ginzburg--Landau equations}, preprint 2004.

\end{thebibliography}
\end{document}